\documentclass[10pt,reqno]{amsart}
\title[Hausdorff measure zero and compactness of the $\mathbf{\dbar}$-Neumann operator]{Hausdorff-$\mathbf{(2n-2)}$ dimensional measure zero set and compactness of the $\dbar$-Neumann operator on $\mathbf{(0,n-1)}$ forms}

\author{YUE ZHANG}
\address{Department of Mathematics, Building 21 Room 404, Zhe Jiang Normal University, Jin Hua, P.R.China, 321004}
\email{yzhangmath@zjnu.edu.cn}
\date{}


\usepackage{amsmath,amsthm,amsfonts,amssymb}
\usepackage[bookmarks=false]{hyperref}
\usepackage{enumerate}

\chardef\bslash=`\\ 

\newtheorem{thm}{Theorem}[section]

\newtheorem{prop}[thm]{Proposition}

\theoremstyle{definition}
\newtheorem{defn}{Definition}[section]

\theoremstyle{remark}
\newtheorem{rem}{Remark}[section]

\newcommand{\thmref}[1]{Theorem~\ref{#1}}

\newcommand{\dbar}{\overline{\partial}}
\newcommand{\dbarad}{\overline{\partial}^*}
\newcommand{\dbaradp}{\overline{\partial}^*_{\varphi}}

\newcommand{\sumprime}{\sideset{}{'}\sum}
\newcommand{\pq}{(P_q)}
\newcommand{\pqt}{(\widetilde{P_q})}

\newcommand{\dom}{\textrm{dom}}

\newcommand{\eval}[2][\right]{\relax
  \ifx#1\right\relax \left.\fi#2#1\rvert}

\begin{document}

\markboth{Hausdorff measure zero and compactness of the $\dbar$-Neumann operator}
{Hausdorff measure zero and compactness of the $\dbar$-Neumann operator}
\begin{abstract}
By using a variant Property $(P_q)$ of Catlin, we discuss the relation of small set of weakly pseudoconvex points on the boundary of pseudoconvex domain and  compactness of the $\dbar$-Neumann operator.  In particular, we show that if the Hausdorff $(2n-2)$-dimensional measure of the weakly pseudoconvex  points on the boundary of a smooth bounded pseudoconvex domain is zero, then the $\dbar$-Neumann operator $N_{n-1}$ is compact on $(0,n-1)$-level $L^2$-integrable forms.
\end{abstract}
\maketitle
\renewcommand{\sectionmark}[1]{}

\section{Introduction}
On a bounded pseudoconvex domain $\Omega$ in $\mathbb{C}^n$, an important question in the $\dbar$-Neumann problem is to study whether there exists a bounded inverse of the complex Laplacian $\square_q=\dbar_{q-1}\dbarad_{q-1}+\dbarad_q\dbar_q$ on the $L^2$-integrable $(0,q)$-type forms of the domain $\Omega$ $(1\leq q\leq n)$ and discuss the regularity property of the inverse if it exists.  To be precise, given a $L^2$-integrable $(0,q)$ form $v$ on $\Omega$, the $\dbar$-Neumann problem is to find $u\in \dom(\square_q)$ such that $\square_q u=v$ and further study regularity property of the solution operator on $L^2$-integrable forms.  We call the (bounded) inverse of $\square_q$ as the $\dbar$-Neumann operator and denote it as $N_q$.  For classical results about the regularity properties of $N_q$, one may check \cite{chenshaw}, \cite{fbkohn}, \cite{hor2}, \cite{stra1} and \cite{zam}.

In this paper, we focus on the study of compactness of the $\dbar$-Neumann operator on specific level forms.  In this regard, Kohn and Nirenberg (\cite{kohnner}) proved that compactness of $N_q$ implies the global regularity of $N_q$ on smooth bounded pseudoconvex domains, here the global regularity means  that $N_q$ maps the space of forms with components smooth up to the boundary of $\Omega$ to itself.  It is well known that compactness of $N_q$ is equivalent to a quantified estimate on $L^2$-integrable forms (see section 2), hence analysis on  compactness of $N_q$ is more robust and has its own interest.  For useful applications of such analysis results, one can check \cite{cat4}, \cite{siqi2}, \cite{hefer}, \cite{henkin}, \cite{sal}, \cite{ven} and references there.

Within the viewpoint of potential analysis theory, there are numerous sufficient conditions for compactness of $N_q$ on a smooth bounded pseudoconvex domain.  For instance, Property $\pq$ in Catlin's work (\cite{cat5}) and Property $(\widetilde{P_q})$ in McNeal's work (\cite{mcneal2}) are well known so far.  In \cite{zhang2}, the author introduced several variant conditions of Property $\pq$ and Property $(\widetilde{P_q})$, which also imply compactness of $N_q$ on high level $L^2$-integrable forms on a smooth bounded pseudoconvex domain.  These variant conditions are obtained by proving a unified estimate  of the twisted Kohn-Morrey-H\"{o}rmander estimate (see in \cite{mcneal2} or section 2.6 in \cite{stra1}) and the $q$-pseudoconvex Ahn-Zampieri estimate (see section 1.9 in \cite{zam} or \cite{ahn}) on a smooth bounded domain.  

In this article, we focus on applying the conditions in \cite{zhang2} on $(0,n-1)$ forms and we discuss the relation of small set of infinite-type points on the boundary of pseudoconvex domain and  compactness of the $\dbar$-Neumann operator $N_{n-1}$.  

This subject  is motivated by the results of Sibony (\cite{sibony}) and Boas (\cite{boas3}) on general pseudoconvex domains: let $q=1$ and assume that the set $K$ of the weakly pseudoconvex points on the boundary $b\Omega$ has Hausdorff $2$-dimensional measure zero in $\mathbb{C}^n$, then the $\dbar$-Neumann operator $N_1$ is compact on $L^2_{(0,1)}(\Omega)$. Boas (\cite{boas3}) has an explicit construction of the function $\lambda$ involved in the proof.  Due to the lack of biholomorphic invariance on Property  $\pq$ when $q>1$, the approach can not be generalized to the case $q>1$ and hence $N_q$ is not known to be compact in the $q>1$ case.

By applying the variant Property $\pq$ when $q=n-1$ in \cite{zhang2}, we prove the following theorem which generalizes above result of Sibony and Boas to the case of $q=n-1$:
\begin{thm}\label{mainthm}
Let $\Omega$ be a smooth bounded pseudoconvex domain in $\mathbb{C}^n$.  If the  Hausdorff $(2n-2)$-dimensional measure of   weakly pseudoconvex points of $b\Omega$ is zero, then the $\dbar$-Neumann operator $N_{n-1}$ is compact on $L^2_{(0,n-1)}(\Omega)$ forms.
\end{thm}
Our result on the $(0,n-1)$-forms is interesting, since under this case, the variant of Property $(P_{n-1})$ we used in proof only involves with the diagonal entries in the complex Hessian, rather than the sum of eigenvalues in the complex Hessian. This fact, in turn, explains why Property $\pq$ of Catlin or Property $(\widetilde{P_q})$ of McNeal is not convenient to apply in the proof of above result.

The paper is organized as follows: in section 2, we list some facts and background materials about the $\dbar$-Neumann problem and related potential analysis results;  in section 3, we prove the main result and mention one example.

\noindent\textbf{Acknowledgment.}  The author wishes to thank Emil Straube and Harold Boas for introducing this problem.

\section{Preliminaries}
Let $L^2_{(0,q)}(\Omega)$ be the space of $(0,q)$-forms ($1 \leq q \leq n$) with $L^2$-integrable coefficients on a bounded domain  $\Omega$ in $\mathbb{C}^n$ ($n\geq 2$).  The $L^2$-norm of a $(0,q)$-form $u$ is defined as $\Vert\sum_{J}^{'}{u_J d\bar{z}_J}\Vert^2=\sum^{'}_{J}\int_{\Omega}{\vert u_J\vert^2 dV(z)}$.  Similarly, the weighted $L^2$-norm of $u$ is defined by $\Vert\sum_{J}^{'}{u_J d\bar{z}_J}\Vert_{\varphi}^2=\sum^{'}_{J}\int_{\Omega}{\vert u_J\vert^2e^{-\varphi} dV(z)}$, where $\varphi\in C^1(\bar{\Omega})$.  Define $\dbar:L^2_{(0,q)}(\Omega)\rightarrow L^2_{(0,q+1)}(\Omega)$ by:
$\dbar(\sumprime_{J}{u_J d\overline{z}_J})=\displaystyle\sum_{j=1}^n\sumprime_J\frac{\partial u_J}{\partial \overline{z}_j}d\overline{z}_j\wedge d\overline{z}_J$.    Let $\textrm{dom}(\dbar)=\lbrace u\in L^2_{(0,q)}(\Omega)|\dbar u\in L^2_{(0,q+1)}(\Omega)\rbrace$ and $\textrm{dom}(\dbarad)=\lbrace v\in L^2_{(0,q+1)}(\Omega)|\exists C>0,~|(v,\dbar u)|\leq C ||u||,\forall u\in \textrm{dom}(\dbar) \rbrace$ be the domain of $\dbar$ and $\dbarad$ respectively.  The weighted $\dbar$-complex is defined similarly in the weighted $L^2$-integrable forms.  We denote the resulting adjoint by $\dbaradp$ and its domain is $\dom(\dbaradp)$.  It is well known that $\dom(\dbaradp)=\dom(\dbarad)$ if $\varphi\in C^1(\bar{\Omega})$. The formal adjoint of $\dbar$ is $\vartheta_\varphi$ such that $(u,\dbar v)_\varphi=(\vartheta_\varphi u,v)_\varphi$ for every $C^\infty$ smooth compactly supported form $v$ on $\Omega$. And $\dbaradp u=\vartheta_\varphi u$ if $u\in\dom(\dbaradp)$.

Given a boundary point $P$ of $\Omega$, we choose vector fields $L_1,\cdots,L_{n-1}$ of type $(1,0)$  which are orthonormal and span $T_z^\mathbb{C}(b\Omega_\epsilon)$ for $z$ near $P$, where $\Omega_\epsilon=\lbrace z\in\Omega|\rho(z)<-\epsilon\rbrace$.  $L_n$ is defined to be the complex normal which can be normalized to be $1$ on the boundary.  We use above vector fields to induce a special boundary chart such that $\lbrace\omega_j\rbrace_{j=1}^n$ is the dual basis of $\lbrace L_j\rbrace_{j=1}^n$ near $P$.  It is then clear that $\dbar f=\sum_{j=1}^n (\bar{L}_jf)\bar{\omega}_j$ for a $C^1$ smooth function $f$.  Let $c_{jk}^i$ be defined by 
$\dbar\omega_i=\sum_{j,k}^n c_{jk}^i\bar{\omega}_j\wedge\omega_k.$
\begin{defn}\label{defn_fjk}
let $f$ be a $C^2$ smooth function, define $f_{jk}=L_j\bar{L}_k f+\sum_i\bar{c}_{jk}^i\bar{L}_i f$. 
\end{defn}
It is then clear that $\partial\dbar f=\sum_{j,k} f_{jk}\omega_j\wedge\bar{\omega}_k$.  
For a general $L^2$-integrable form $u=\sumprime_{|J|=q}u_J\bar{\omega}_J$ in the special boundary chart, we have:
\begin{equation}\label{eqn1}
\dbar u=\sumprime_{|K|=q-1}\sum_{i<j}\left( \bar{L}_i u_{jK}-\bar{L}_j u_{iK} \right) \bar{\omega}_i\wedge\bar{\omega}_j\wedge\bar{\omega}_K+\cdots,
\end{equation}

\begin{equation}\label{eqn2}
\vartheta_\varphi u=-\sumprime_{|K|=q-1}\sum_{j\leq n}\delta_{\omega_j}(u_{jK})\bar{\omega}_K+\cdots,
\end{equation}
where $\delta_{\omega_j} f=e^\varphi L_j(e^{-\varphi}f)$ for $L^2$-integrable functions $f$.  The dots in above two equations are the terms that only involve with the coefficients of $u$ and the differentiation of the coefficients of $L_j$ or $\bar{\omega}_K$.

We define the complex Laplacian as $\square_q u:=\dbarad\dbar u+\dbar\dbarad u$ on $L^2_{(0,q)}$ forms.  Here we suppress the subscript of the level of the form in $\dbar$ and $\dbarad$ for simplicity.  We call the  inverse operator of $\square_q$ as the $\dbar$-Neumann operator, and denote it as $N_q$.  H{\"o}rmander (\cite{hor1,hor2})showed that $\square_q$ has a bounded inverse $N_q$ on $L^2_{(0,q)}(\Omega)$ when $\Omega$ is a bounded pseudoconvex domain.  $N_q$ is said to be compact on $L^2_{(0,q)}(\Omega)$ if the image of the unit ball in $L^2_{(0,q)}(\Omega)$ under $N_q$ is relatively compact in $L^2_{(0,q)}(\Omega)$.  We can characterize the compactness of $N_q$ by the following well known fact (see \cite{mcneal2} or \cite{stra1}, Proposition 4.2):

\begin{prop}\label{new_section1_1}
Let $\Omega$ be a bounded pseudoconvex domain in $\mathbb{C}^n$, $1\leq q\leq n$.  Then the following are equivalent:
\begin{enumerate}[{\normalfont(i)}]
\item $N_q$ is compact as an operator on $L^2_{(0,q)}(\Omega)$.
\item For every $\epsilon>0$, there exists a constant $C_\epsilon$ such that we have the compactness estimate:
\[||u||^2\leq\epsilon(||\dbar u||^2+||\dbarad u||^2)+C_\epsilon||u||^2_{-1}~\textrm{for}~u\in \dom(\dbar)\cap \dom(\dbarad).\]
\item The canonical solution operators $\dbarad N_q:L^2_{(0,q)}(\Omega)\cap\textrm{ker}(\dbar)\rightarrow L^2_{(0,q-1)}(\Omega)$ and $\dbarad N_{q+1}:L^2_{(0,q+1)}(\Omega)\cap\textrm{ker}(\dbar)\rightarrow L^2_{(0,q)}(\Omega)$ are compact.
\end{enumerate}
\end{prop}

Catlin (\cite{cat5}) showed that if $\Omega$ be a smooth bounded pseudoconvex domain and $b\Omega$ satisfies Property $(P_q)$, then $N_q$ is compact on $L^2_{(0,q)}(\Omega)$.  McNeal (\cite{mcneal2}) showed that Property $(P_q)$ can be weakened to Property $\pqt$ on individual function level, and still implies compactness of $N_q$.  We list the definition of Property $(P_q)$ here for use in section 3:
\begin{defn}\label{defn_pq}
A compact set $K\subset\mathbb{C}^n$  has Property $(P_q)$ ($1\leq q\leq n$) if for any $M>0$, there exists an open neighborhood $U$ of $K$ and a $C^2$ smooth function $\lambda$ on $U$ such that $0\leq\lambda\leq 1$ on $U$ and $\forall z\in U$, the sum of any $q$ eigenvalues of the complex Hessian $\left(\frac{\partial^2\lambda}{\partial z_j\partial \bar{z}_k}\right)_{j,k}$ is at least $M$.
\end{defn}

In \cite{zhang2}, the author introduced several variant conditions of Property $(P_q)$ and Property $\pqt$ which still imply compactness of $N_q$ on smooth bounded pseudoconvex domains.  We list the definition of  a variant of Property $(P_{n-1})$ in \cite{zhang2}, which will be used in this article.
\begin{defn}
For a smooth bounded pseudoconvex domain $\Omega\subset\mathbb{C}^n$ $(n>2)$, $b\Omega$ has Property $(P_{n-1}^\#)$ if there exists a finite cover $\lbrace V_j\rbrace_{j=1}^N$ of $b\Omega$ with special boundary charts  and the following holds on each $V_j$:  for any $M>0$, there exists a neighborhood $U$ of $b\Omega$ and a $C^2$ smooth function $\lambda$ on $U\cap V_j$, such that $0\leq\lambda(z)\leq 1$ and there exists $t$ ($1\leq t\leq n-1$) such that $\lambda_{tt}\geq M$ on $U\cap V_j$.
\end{defn}
Here, as in Definition \ref{defn_fjk}, $\lambda_{tt}=L_t\bar{L}_t \lambda+\sum_i\bar{c}_{tt}^i\bar{L}_i \lambda$ is the diagonal entry in the Hessian matrix $(\lambda_{jk})$.  We have the following result in \cite{zhang2}:
\begin{thm}[\cite{zhang2}]\label{zhangthm1}
Let $\Omega\subset\mathbb{C}^n$ $(n>2)$ be a smooth bounded pseudoconvex domain.  If $b\Omega$ has Property $(P_{n-1}^\#)$, then the $\dbar$-Neumann operator $N_{n-1}$ is compact on $L^2_{(0,n-1)}(\Omega)$.
\end{thm}
We  also need the following result due to Sibony (\cite{sibony}):
\begin{prop}\label{sibony2}
Let $K$ be a compact subset in $\mathbb{C}^n$ ($n\geq 1$) and $K$ has Lebesgue measure zero in $\mathbb{C}^n$.  Then $K$ has Property $(P_n)$ in $\mathbb{C}^n$.
\end{prop}
The original result is formulated for $n=1$ case. But the sum of any $n$ eigenvalues of the complex Hessian of $\lambda$ in $\mathbb{C}^n$ is equal to the real Laplacian of $\lambda$ in $\mathbb{R}^{2n}$, and most of the classical potential results which were used in the proof of this result can also be  formulated in $\mathbb{R}^{2n}$, hence the result can be generalized to $n>1$ case trivially.

\section{Proof of main theorem}

\begin{proof}[Proof of \thmref{mainthm}]

Let $\lbrace \xi_j\rbrace_{j=1}^{n-1}$ be the orthonormal coordinates which span the complex tangent space $Z$ in the special boundary chart at a boundary point $P$.  Let $V$ be a neighborhood of the boundary point $P$, and $K$ be the weakly pseudoconvex points on the boundary $b\Omega$.  Let $\pi^{Z}:\mathbb{C}^n\to\mathbb{C}^{n-1}$ be the projection map from $\mathbb{C}^n$ onto the complex tangent space $Z$ at $P$.

The set $\pi^{Z}(K\cap V)$ has Hausdorff-$(2n-2)$ dimensional measure zero in a copy of $\mathbb{C}^{n-1}$, since any continuous map preserves Hausdorff measure zero set.  Since Hausdorff-$(2n-2)$ dimensional measure is equivalent to Lebesgue measure in $\mathbb{C}^n$ (modulo a constant), by Proposition \ref{sibony2}, the set $\pi^{Z}(K\cap V)$ has Property $(P_{n-1})$ of Catlin.  That is, for any $M>0$, there exists a neighborhood in $\mathbb{C}^{n-1}$ of $\pi^{Z}(K\cap V)$ and a $C^2$ smooth function $\lambda^M(\xi_1,\cdots,\xi_{n-1})$ such that $0\leq \lambda^M\leq 1$ and the real Laplacian $\Delta \lambda^M(\xi_1,\cdots,\xi_{n-1})\geq M$ on the above neighborhood of $\pi^{Z}(K\cap V)$.  Here the Laplacian is taken with respect to the coordinates $(\xi_1,\cdots,\xi_{n-1})$ in $\mathbb{C}^{n-1}$.  Define $\lambda^M_{jk}$ which is same in Definition \ref{defn_fjk}, therefore $\Delta \lambda^M(\xi_1,\cdots,\xi_{n-1})=\sum_{j=1}^{n-1}\lambda^M_{jj}$ by using the invariance of real Laplacian under orthonormal coordinates change.  


On the neighborhood $V$, define the trivial extension function $\eta^M(\xi_1,\xi_2,\cdots,\xi_n)=\lambda^M(\xi_1,\cdots,\xi_{n-1})$. Then the real Laplacian $\Delta\eta^M$ on the boundary is equal to the real Laplacian $\Delta\lambda^M$.  Consider the entries in the complex Hessian of $(\eta^M_{jk})$, the size of this matrix is $n\times n$.  For $1\leq j\leq n-1$, $\eta^M_{jj}=\lambda^M_{jj}$ by using Definition \ref{defn_fjk}.

Now let the set $E^M_{j}=\pi^{Z}(K\cap V) \cap \lbrace \eta^M_{jj}\geq \frac{M}{n-1} \rbrace$, $1\leq j\leq n$.  By definition of $\lambda^M$, we have $\pi^{Z}(K\cap V)\subseteq \bigcup_{j=1}^{n-1}E_j^M$.  Then $\bigcup_{j=1}^{n-1} \big( \pi^{-1}_Z(E^M_j) \cap V\big)\supseteq K\cap V$, here $\pi^{-1}_Z$ is the inverse map of $\pi^{Z}$.


The diagonal entry $\eta^M_{jj}$ in the complex Hessian of $(\eta^M_{j,k})$ satisfies the conditions in the definition of Property $(P_{n-1}^\#)$ on each $\pi^{-1}_Z(E^M_j) \cap V$ when $1\leq j\leq n-1$.  Now since $\bigcup_{j=1}^{n-1} \big( \pi^{-1}_Z(E^M_j) \cap V\big)\supseteq K\cap V$ by the previous paragraph, we can apply Property $(P_{n-1}^\#)$ together with partition of unity to prove the compactness estimate locally on $V$.  The cut-off functions in the partition should produce extra partial derivatives by hitting $\dbar$ and $\dbar$, but those derivatives can be handled in the same way as the proof of \thmref{zhangthm1}, hence the desired compactness estimate (see (ii) in Proposition \ref{new_section1_1}) will not be affected.  Also for the strongly pseudoconvex points on $V$, they are naturally of D'Angelo's finite type and hence compactness estimate holds there (see \cite{cat5}, \cite{chenshaw} or \cite{stra1}).  Since compactness of the $\dbar$-Neumann operator is a local property, the conclusion follows.
\end{proof}

\begin{rem}
For the case of Hausdorff $2$-dimensional measure and compactness of $N_1$, as we pointed out in the introduction section, the essential argument in Sibony and Boas's work (\cite{boas3} and \cite{sibony}) is to show that the infinite-type points on the boundary satisfy Property $(P_1)$.  In such argument, the idea is to project the set $K$ of infinite-type points to each $z_j$-plane and the resulting set satisfies Property $(P_1)$ on each complex $1$-dimensional plane, hence summing all involved functions in the definition of Property $(P_1)$ will give the desired conclusion.  Now in our case of \thmref{mainthm}, such summation of functions does not work since eigenvalues from each respective complex Hessian interfere the summation of eigenvalues in the whole complex Hessian.  Therefore, verifying Property $(P_q)$ or Property $\pqt$ under such case appears not to work.  A detailed explanation of such phenomenon under potential analysis background can also be found in the author's recent work (see remarks after Corollary 3.2 in \cite{zhang1}).  

Our result in \thmref{mainthm} shows that small set of weakly pseudoconvex points (or infinite-type points) on the boundary in the sense of Hausdorff-$(2n-2)$ dimensional measure is benign in the compactness of $N_{n-1}$.  When $1<q<n-1$, whether similar conclusion holds in the sense of Hausdorff-$2q$  dimensional measure is not known yet.  In such case, a certain arrangement on projections onto each $q$-dimensional subspace needs to be found.

\end{rem}
For an example when \thmref{mainthm} holds, we give one example from \cite{zhang1} and refer the reader to there for details of calculation.
\begin{prop}
Define a smooth complete Hartogs domain $\Omega\subset\mathbb{C}^3$ by:
$$\Omega=\lbrace(z_1,z_2,z_3)|~|z_3|^2<e^{-\varphi(z_1)-\psi(z_2)},z_1\in\mathbb{D}(0,1),z_2\in\mathbb{D}(0,1)\rbrace.$$   Assume that $\varphi,\psi\in C^\infty(\mathbb{D}(0,1))$ and subharmonic on $\mathbb{D}(0,1)$ in the respective complex plane.  Assume further that the boundary points $(z_1,z_2,z_3)$ are strictly pseudoconvex when $(z_1,z_2)$ is close to $b(\mathbb{D}(0,1)\times\mathbb{D}(0,1))$.  If the Hausdorff $4$-dimensional measure of the weakly pseudoconvex points of $b\Omega$ is zero, then the $\dbar$-Neumann operator $N_2$ is compact. 
\end{prop}


\begin{thebibliography}{100}
\bibitem{ahn} H. Ahn, Global boundary regularity for the $\dbar$-equation on $q$-pseudoconvex domains.  \emph{Math. Nachr.} 280 (2007), 343-350.

\bibitem{boas3} H.\,P. Boas, Small sets of infinite type are benign for the $\dbar$-Neumann problem. \emph{Proc. Amer. Math. Soc.} \textbf{103}, 569-578 (1988).


\bibitem{cat5} D. Catlin, Global regularity of the $\dbar$-Neumann problem.  In \emph{Complex analysis of several variables} (Madison, 1982).  Proc.  Sympos. Pure Math. 41, Amer. Math. Soc., Providence 1984, 39-49.


\bibitem{cat4} D. W. Catlin and J. P. D'Angelo, Positivity conditions for bihomogeneous polynomials.  \emph{Math. Res. Lett.} 4 (1997), 555-567.

\bibitem{chenshaw} S.-C. Chen and M.-C. Shaw, \emph{Partial differential equations in several complex variables}.  AMS/IP Stud. Adv. Math. 19, Amer. Math. Soc., Providence, R.I., 2001.



\bibitem{fbkohn} G. B. Folland and J. J. Kohn, \emph{The Neumann problem for the Cauchy-Riemann complex}.  Ann. of Math. Stud. 75, Princeton University Press, Princeton, N.J., 1972.


\bibitem{siqi2} S. Fu and E. J. Straube, Compactness in the $\dbar$-Neumann problem.  In \emph{Complex analysis and geometry} (Columbus 1999).  Ohio State Univ. Math. Res. Inst. Publ. 9, Walter de Gruyter, Berlin 2001, 141-160.




\bibitem{hefer} T. Hefer and I. Lieb, On the compactness of the $\dbar$-Neumann operator.  \emph{Ann. Fac. Sci. Toulouse Math. (6)} 9 (2000), 415-432.



\bibitem{henkin} G. M. Henkin and A. Iordan, Compactness of the Neumann operator for hyperconvex domains with non-smooth $B$-regular boundary.  \emph{Math. Ann.} 307 (1997), 151-168.

\bibitem{hor1} L. H{\"o}rmander, $L^2$ estimates and existence theorem for the $\dbar$ operator. \emph{Acta Math}. 113 (1965), 89-152.

\bibitem{hor2} L. H{\"o}rmander, \emph{An introduction to complex analysis in several variables}. 3rd ed., North-Holland Math. Library 7, North-Holland Publishing Co., Amsterdam 1990.







 
\bibitem{kohnner} J. J. Kohn and L. Nirenberg, Non-coercive boundary value problems.  \emph{Comm. Pure Appl. Math.} 18 (1965), 443-492.



\bibitem{mcneal2} J. McNeal, A sufficient condition for compactness of the $\dbar$-Neumann operator.  \emph{J. Funct. Anal.} 195 (2002), 190-205.




\bibitem{sal} N. Salinas, Noncompactness of the $\dbar$-Neumann problem and Toeplitz $C^*$-algebras.  In \emph{Several complex variables and complex geometry}, Part 3 (Santa Cruz, CA, 1989), Proc. Sympos. Pure Math. 52, Amer. Math. Soc., Providence, R.I., 1991, 329-334.

\bibitem{sibony} N. Sibony, Une classe de domaines pseudoconvexes.  \emph{Duke Math. J.} 55 (1987), 299-319.


\bibitem{stra1} E. J. Straube, \emph{Lectures on the $L^2$-Sobolev Theory of the $\dbar$-Neumann Problem}.  ESI Lectures in Mathematics and Physics, European Math. Society Publishing House, Z{\"{u}}rich, 2010.

\bibitem{zam}  G. Zampieri, \emph{Complex analysis and CR geometry}.  University  Lecture Ser. 43, Amer. Math. Soc., Providence, R.I., 2008.

\bibitem{ven} U. Venugopalkrishna, Fredholm operators associated with strongly pseudoconvex domains in $\mathbb{C}^n$. \emph{J. Funct. Anal.} 9 (1972), 349-373.

\bibitem{zhang1} Y. Zhang, Some aspects of Property $\pq$ for $q>1$, \emph{Math. Nachr.}, 290 (2017), 1119-1134.

\bibitem{zhang2} Y. Zhang, {A sufficient condition for the compactness of the $\dbar$-Neumann operator on high level forms}, submitted. {\tt arXiv:1710.09614 [math.CV]}

\end{thebibliography}
\end{document}